\documentclass[11pt]{amsart}

\usepackage{epigamath}

\usepackage[english]{babel}

\numberwithin{equation}{section}

\usepackage[shortlabels]{enumitem}
\setlist[enumerate,1]{label={\rm(\roman*)}, ref={\rm\roman*}}

\usepackage{amscd, amsmath, amssymb, fancyhdr, graphicx, epsfig, url, caption}
\usepackage{cleveref}

\usepackage{url}

\newtheorem{theorem}{Theorem}[section]
\newtheorem{corollary}[theorem]{Corollary}
\newtheorem{lemma}[theorem]{Lemma}
\newtheorem{proposition}[theorem]{Proposition}
\newtheorem{conjecture}[theorem]{Conjecture}

\theoremstyle{definition}
\newtheorem{definition}[theorem]{Definition}

\theoremstyle{remark}

\newtheorem{remark}[theorem]{Remark}

\newcommand{\Z}{{\mathbb Z}}
\newcommand{\C}{{\mathbb C}}

\newcommand{\R}{{\mathbb R}}
\newcommand{\Q}{{\mathbb Q}}
\def\1{\sqrt{-1}\:}

\renewcommand{\bar}{\overline}
\renewcommand{\epsilon}{\varepsilon}
\renewcommand{\leq}{\leqslant}

\newcommand{\Teich}{\mathcal{T}}
\newcommand{\nef}{\mathrm{nef}}
\newcommand{\sa}{\mathrm{sa}}

\newcommand{\MCG}{\mathcal{G}}
\newcommand{\MBM}{\mathrm{MBM}}
\newcommand{\Comp}{\mathscr{I}}
\newcommand{\Per}{\mathcal{P}er}
\newcommand{\Perspace}{\mathcal{D}}
\newcommand{\im}{\operatorname{im}}
\newcommand{\Pic}{\operatorname{Pic}}

\newcommand{\Diff}{\operatorname{\mathscr Diff}}

\newcommand{\Tw}{\operatorname{Tw}}
\renewcommand{\Re}{\operatorname{Re}}

\newcommand{\bbA}{\mathbb{A}}
\newcommand{\bbZ}{\mathbb{Z}}
\newcommand{\bbQ}{\mathbb{Q}}
\newcommand{\bbR}{\mathbb{R}}
\newcommand{\bbC}{\mathbb{C}}
\newcommand{\bbP}{\mathbb{P}}
\newcommand{\CC}{\mathcal{C}}
\newcommand{\DD}{\mathcal{D}}
\newcommand{\HH}{\mathcal{H}}
\newcommand{\KK}{\mathcal{K}}
\newcommand{\LL}{\mathcal{L}}
\newcommand{\MM}{\mathcal{M}}
\newcommand{\OO}{\mathcal{O}}
\newcommand{\TT}{\mathcal{T}}
\newcommand{\VV}{\mathcal{V}}
\renewcommand{\ge}{\geqslant}
\renewcommand{\le}{\leqslant}

\newcommand{\st}{\mid} 

\newcommand{\cnv}{\,\lrcorner\,}
\newcommand{\wdg}{\wedge}

\newcommand{\lra}{\longrightarrow}
\def\lowsim{\vbox to 0pt{\vss\hbox{$\scriptstyle\sim$}\vskip-1.6pt}}

\EpigaVolumeYear{9}{2025} \EpigaArticleNr{26} \ReceivedOn{September 30, 2024}
\InFinalFormOn{April 1, 2025}
\AcceptedOn{July 13, 2025}

\title{The abundance and SYZ conjectures in families of hyperk\"ahler manifolds}
\titlemark{Abundance and SYZ in families of HK manifolds}

\author{Andrey Soldatenkov}
\address{Universidade Estadual de Campinas,
Departamento de Matem\'atica - IMECC,
Rua S\'ergio Buarque de Holanda, 651,
CEP 13083-859, Campinas, SP, Brazil}
\email{aosoldatenkov@gmail.com}

\author{Misha Verbitsky}
\address{Instituto Nacional de Matem\'atica Pura e
Aplicada (IMPA), Estrada Dona Castorina, 110,
Jardim Bot\^anico, CEP 22460-320,
Rio de Janeiro, RJ - Brazil}
\email{verbit@impa.br}

\authormark{A.~Soldatenkov and M.~Verbitsky}

\AbstractInEnglish{Let $L$ be a holomorphic line bundle on a hyperk\"ahler manifold $M$, with $c_1(L)$ nef and not big. The SYZ conjecture predicts that $L$ is semiample. We prove that this is true, assuming that $(M,L)$ has a deformation $(M',L')$ with $L'$ semiample.  We introduce a version of the Teichm\"uller space that parametrizes pairs $(M,L)$ up to isotopy. We prove a version of the global Torelli theorem for such Teichm\"uller spaces and use it to deduce the deformation invariance of semiampleness.}



\acknowledgement{Misha Verbitsky was partially supported by 
FAPERJ SEI-260003/000410/2023 and CNPq - Process 310952/2021-2.\\
Andrey Soldatenkov was partially supported by 
FAPESP grant 2024/23819-0.}

\begin{document}



\maketitle

\begin{prelims}

\DisplayAbstractInEnglish

\bigskip

\DisplayKeyWords

\medskip

\DisplayMSCclass

\end{prelims}


\newpage

\setcounter{tocdepth}{1}

\tableofcontents


\section{Introduction}

The present paper has its roots in two conjectures: the abundance conjecture with its generalizations, coming from birational geometry, and the SYZ conjecture, which comes from physics, the theory of calibrations and Calabi--Yau geometry. Let us briefly recall these conjectures.

\subsection{The abundance conjecture}

Recall that a line bundle $L$ on a compact K\"ahler manifold $M$ is called {\em nef}\, if $c_1(L)$ belongs to the closure of the K\"ahler cone of $M$, and {\em big} if its Iitaka dimension $\kappa(L)$ (see Section~\ref{sec_semiample} for a definition) equals the dimension of $M$. The bundle $L$ is called {\em semiample} if there exists a $k > 0$ such that $L^k$ is generated by its global sections.

Semiampleness of a line bundle implies that some power of this bundle defines a morphism from $M$ to some other projective variety, and constructing such morphisms is often crucial for the study of the geometry of $M$. It is therefore one of the central problems in algebraic geometry to find sufficient conditions that guarantee the semiampleness of a given line bundle. The following is one of the foundational results in that direction, which we state in a simplified form to streamline the exposition.

\begin{theorem}[Kawamata's base-point-free theorem, \textit{cf.} \protect{\cite[Theorem 6.1]{_Kawamata:Pluricanonical_}}]
Let $L$ be a nef line bundle on a projective manifold $M$ such that $L^{\otimes m}\otimes K_M^{-1}$ is nef and big for some $m>0$, where $K_M$ is the canonical bundle of\,~$M$. Then $L$ is semiample.
\end{theorem}

The above theorem serves as a motivation for a number of conjectures, including the abundance conjecture and its generalizations (see \textit{e.g.}~\cite{_Demailly_Peternell_Schneider:ps-eff_}, Conjecture 2.7.2). One of its versions, again in a simplified form, is stated as follows.

\begin{conjecture}[Generalized abundance conjecture,   \textit{cf.}
\cite{_Lazic_Peternell:gen_abund_I_}]\label{conj_abundance} 
Let $M$ be a projective manifold with pseudoeffective canonical bundle $K_M$.
Let $L$ be a nef line bundle such that $L \otimes K_M^{-1}$ is also nef.
Then $L$ is numerically equivalent to a semiample line bundle $L'$.
\end{conjecture}

In the present paper we will assume that $M$ is a compact hyperk\"ahler manifold (see Section~\ref{_HK_Subsection_} for definitions); therefore, the canonical bundle $K_M$ will be trivial. Kawamata's theorem in this case implies that a big and nef line bundle $L$ is semiample. On the other hand, if we drop the assumption of bigness and assume only that $L$ is nef, then the abundance conjecture stated above still predicts the semiampleness of $L$ (a numerical equivalence class on a hyperk\"ahler manifold consists of a single line bundle because $b_1(M) = 0$). If true, the conjecture implies that a power of $L$ defines a morphism onto a projective variety of dimension smaller than $\dim(M)$, which in fact turns out to be {\em a Lagrangian fibration} (see Section~\ref{sec_semiample} for details). This gives a link to another important conjecture in hyperk\"ahler geometry which we discuss next.

\subsection{The hyperk\"ahler SYZ conjecture}

The hyperk\"ahler SYZ conjecture was formulated many times independently since the 1990s; for its history and its relevance to string physics, see \cite{_Sawon_,_Verbitsky:SYZ_}. Among the first people who stated this conjecture are Tyurin, Bogomolov, Hasset, Tschinkel, Sawon and Huybrechts (\cite{_Hassett_Tschinkel:SYZ_conj_,_Sawon_,_Huybrechts:basic_}).  In its weakest form, it can be stated as follows. 

\begin{conjecture}[The hyperk\"ahler SYZ conjecture, weak form]\label{_SYZ_weak_Conjecture_}
Let $M$ be a hyperk\"ahler manifold. Then $M$ can be deformed to a
hyperk\"ahler manifold admitting a holomorphic Lagrangian fibration.
\end{conjecture}

A holomorphic Lagrangian fibration mentioned in the conjecture is a morphism with connected fibres $\pi\colon M \to B$ to a normal projective variety $B$ with $\dim(B) = \frac{1}{2}\dim(M)$ and all fibres of $\pi$ being Lagrangian subvarieties of $M$ (see Section~\ref{sec_semiample}).  Denote by $q$ the BBF form on $H^2(M, \bbQ)$ (see Section~\ref{_HK_Subsection_}).  If $\OO_B(1)$ is an ample line bundle on $B$, then $L = \pi^*\OO_B(1)$ is semiample and satisfies $q(c_1(L))=0$, the latter assertion following from the Fujiki relations (\ref{eqn_Fujiki}). Taking into account these observations and~\Cref{conj_abundance}, one arrives at the following more precise version of the hyperk\"ahler SYZ conjecture; see \cite[Conjecture 1.7]{_Verbitsky:SYZ_}. 

\begin{conjecture}[The hyperk\"ahler SYZ conjecture, strong form]\label{_SYZ_strong_}
Let $M$ be a hyperk\"ahler manifold and $L$ a non-trivial nef line bundle on $M$ with $q(c_1(L))=0$. Then there exist a Lagrangian fibration $\pi\colon M \to B$ and an ample line bundle $\OO_B(1)$ on $B$ such that $L^k = \pi^* \OO_B(1)$ for some integer $k >0$.
\end{conjecture}

The hyperk\"ahler SYZ conjecture has been extensively studied, and there exists substantial evidence that the conjecture should be true; let us mention in particular the results of \cite{_COP:non-alge_}, where the conjecture is proven for non-algebraic hyperk\"ahler manifolds of algebraic dimension $\frac 1 2 \dim_\C M$.  The paper \cite{AH} gives an overview of other partial results in this direction.

In each of the currently known deformation classes of hyperk\"ahler manifolds, there exist manifolds admitting Lagrangian fibrations.  Therefore, a natural approach to~\Cref{_SYZ_strong_} (at least for the known deformation classes) is to obtain the result by deformation techniques. This approach was extensively explored by Matsushita; see \cite{_Matsushita:nef_,_Matsushita:families_,_Matsushita:divisors_}.  Using a theorem of Voisin \cite{_Voisin:Lagrang_def_} about deformations of Lagrangian submanifolds of a hyperk\"ahler manifold $M$, Matsushita showed in \cite{_Matsushita:families_} that if a line bundle $L$ with $q(c_1(L))=0$ is semiample, then it remains semiample after a small deformation of the pair $(M, L)$; \textit{i.e.} semiampleness is open in families.

In \cite[Theorem 1.2]{_Matsushita:divisors_} Matsushita extends his results and obtains more precise statements about the behaviour of semiampleness under local deformations, assuming additionally that the base of the Lagrangian fibration defined by the semiample line bundle is isomorphic to $\bbC P^n$. In particular, as part of his proof, see \cite[Claim 3.2]{_Matsushita:divisors_}, Matsushita obtains deformation invariance of semiampleness under the assumption that the base is isomorphic to $\bbC P^n$.

In the literature one often finds discussions of a weaker ``birational version'' of~\Cref{_SYZ_strong_}. In this birational version one does not assume the bundle $L$ to be nef (it should only be in the closure of the birational K\"ahler cone), and the conclusion should be that $M$ is bimeromorphic to another hyperk\"ahler manifold $M'$ on which~$L$ induces a Lagrangian fibration. The works of Matsushita combined with other results imply a positive solution of the birational SYZ conjecture for the known deformation types of hyperk\"ahler manifolds.

\begin{remark}\label{_SYZ_examples_ref_Remark_}
  Here are the references to the proof of the bimeromorphic version of SYZ: for the Hilbert schemes of points on K3 surfaces \cite[Theorem 1.5]{_Bayer_Macri_} and \cite[Theorems 1.3 and~6.3]{_Markman:SYZ_}, for the deformations of the generalized Kummer varieties \cite[Proposition 3.38]{_Yoshioka_}, for the O'Grady's sixfold \cite[Corollary 1.3 and~7.3]{_Mongardi_Rapagnetta_}), and for the O'Grady's tenfold \cite[Theorem 2.2]{_Mongardi_Onorati_}.
\end{remark}

\subsection{Results}

The main goal of this paper is to complete the study of the behaviour of semiampleness under deformation initiated by Matsushita. Recall that a line bundle $L$ is called {\em isotropic} if $q(c_1(L))=0$, where $q$ is the BBF quadratic form (Section~\ref{_HK_Subsection_}).  We build the global moduli theory of complex structures of hyperk\"ahler type equipped with isotropic semiample line bundles.  This leads to the positive answer to~\Cref{_SYZ_strong_} under the assumption that the pair $(M,L)$ admits at least one deformation $(M',L')$ with semiample $L'$. This is our~\Cref{thm_abundance}.

In order to prove this result, we develop a moduli theory for the pairs $(M,L)$ where $M$ is a hyperk\"ahler manifold and $L$ a semiample line bundle of vanishing BBF square on $M$.  Fixing a cohomology class $\ell \in H^2(M, \Z)$, we define a semiample Teichm\"uller space $\Teich_{\sa}(M,\ell)$ parametrizing pairs $(M, L)$ where $L$ is semiample with $c_1(L) = \ell$.  The space $\Teich_{\sa}(M,\ell)$ is defined as a subset of $\Teich(M,\ell)\subset \Teich(M)$, where $\Teich(M,\ell)$ is a divisor in the usual Teichm\"uller space where the class $\ell$ stays of Hodge type $(1,1)$; see Section~\ref{sec_Teich_sa} for details.  We define the period domain for $\Teich_{\sa}(M,\ell)$ and study the fibres of the corresponding period map.  Our main result is a version of the global Torelli theorem for $\Teich_{\sa}(M,\ell)$,~\Cref{thm_Torelli_sa}, claiming in particular that the period map is surjective when the semiample Teichm\"uller space is non-empty.

Our results rely only on the openness of semiampleness proved by Matsushita in \cite{_Matsushita:families_} and do not use \cite{_Matsushita:divisors_}. In particular, we do not assume anything about the base of our Lagrangian fibrations, which makes our results independent of the well-known open problems regarding the smoothness of the base. Our main tool is the construction of special families of Lagrangian fibrations called the degenerate twistor families. We rely on our previous result \cite{_SV:Herm_} claiming that all fibres of such families are hyperk\"ahler. This implies that the semiample Teichm\"uller space is covered by affine lines, which is a key step in our proof of~\Cref{thm_Torelli_sa}.

Regarding the application of our results to the known deformation types of hyperk\"ahler manifolds, let us mention the following.  It is shown in \cite[Corollary 3.10]{_KV:Lagra_families_} that there are only finitely many classes of divisors $\Teich (M,\ell)$ up to the action of the mapping class group.  For each of the known deformation types, the equivalence classes of the divisors $\Teich (M,\ell)$ are classified (see \Cref{_SYZ_examples_ref_Remark_}), and in each of these divisors, there is a member for which $\ell=c_1(L)$, where $L$ is a semiample line bundle.~\Cref{thm_abundance} therefore implies~\Cref{_SYZ_strong_} for all known deformation classes of hyperk\"ahler manifolds.

\subsection*{Acknowledgments}
We are grateful to Ekaterina Amerik for her interest and discussions,
and to Dominique Mattei for the useful remarks and comments.


\section{Hyperk\"ahler manifolds and Teichm\"uller spaces}
\label{_hk_intro_Section_}


In this section we recall the necessary notions and known results from hyperk\"ahler geometry.  Then we introduce the semiample Teichm\"uller space, which will be the central object of our study.

\subsection{Basic definitions and conventions}
\label{_HK_Subsection_}

Our setup will be the same as in \cite{_SSV:rigid_currents_,_SV:Moser_,_SV:Herm_}.  For a detailed treatment of the subject, we refer to \cite{_Beauville_,_Huybrechts:basic_,_Verbitsky:Torelli_}.

Let $M$ be a compact simply connected $C^\infty$-manifold and $g$ a Riemannian metric on it. We will say that $g$ is {\em hyperk\"ahler of maximal holonomy} if the holonomy group of the Levi--Civita connection $\nabla^g$ is isomorphic to $\operatorname{Sp}(n)$. In this case the connection $\nabla^g$ preserves a triple of complex structures $I$, $J$, $K$ such that $IJ = -JI = K$, and the tuple $(g,I,J,K)$ is called {\em a hyperk\"ahler structure} on $M$. From now on we will assume that $M$ is hyperk\"ahler of maximal holonomy.

A complex structure $I$ on a compact manifold $M$ is {\em of hyperk\"ahler type} if it is part of a hyperk\"ahler structure for some hyperk\"ahler metric of maximal holonomy. As follows from the Calabi--Yau theorem (see \cite{_Beauville_,_Besse:Einst_Manifo_}), $I$ has this property if and only if it is holomorphically symplectic, of K\"ahler type, and $\pi_1(M)=1$, $H^{2,0}(M)=\C$.  We will denote by $\Comp(M)$ the set of complex structures of hyperk\"ahler type with its natural Fr\'echet topology and by $\Diff^\circ(M)$ the identity connected component of the diffeomorphism group of $M$.  The quotient $\Teich(M)=\Comp(M)/\Diff^\circ(M)$ is called {\em the Teichm\"uller space} of $M$. The quotient group $\MCG(M) = \Diff(M)/\Diff^\circ(M)$ is called {\em the mapping class group} of $M$.  The Teichm\"uller space is a non-Hausdorff complex manifold of dimension $b_2(M) - 2$, and there is a natural $\MCG(M)$-action on $\Teich(M)$.

For a complex structure $I\in \Comp(M)$, the complex manifold $(M, I)$ admits a unique (up to multiplication by a constant) holomorphic symplectic form $\sigma_I\in \Lambda^{2,0}_I M$. It can be written as $\sigma_I=\omega_J+\1\omega_K$, where $\omega_J$ and $\omega_K$ are the K\"ahler forms on $(M,J)$ and $(M,K)$.  A symplectic manifold is even-dimensional, and we define $2n = \dim_\bbC M$.

For a hyperk\"ahler manifold $M$, the group $V_\bbZ = H^2(M,\bbZ)$ is torsion-free, and we define $V_\bbQ = V_\bbZ\otimes \bbQ$, $V_\bbR = V_\bbZ\otimes \bbR$ and $V_\bbC = V_\bbZ\otimes \bbC$. The vector space $V_\bbQ$ carries a quadratic form $q\in S^2V^*_\bbQ$ of signature $(3,b_2(M) - 3)$ called the Beauville--Bogomolov--Fujiki form or {\em the BBF form}. We normalize $q$ so that its restriction to $V_\bbZ$ is integral and primitive.  We will use the same letter $q$ to denote the quadratic form $q\in S^2V^*_\bbQ$ and the bilinear symmetric form $V_\Q \otimes V_\Q \to \Q$.  The BBF form satisfies {\em the Fujiki relations:} there exists a positive constant $c_M\in \bbQ$ such that for any $x\in V_\bbC$
\begin{equation}\label{eqn_Fujiki}
\int_M x^{2n} = c_M q(x)^n.
\end{equation}

It follows from the Fujiki relations that the cohomology class $x = [\sigma_I]$ of the symplectic form satisfies $q(x) = 0$ and $q(x,\bar{x}) > 0$.  Define {\em the period domain}
\begin{equation}\label{eqn_perspace}
\Perspace = \{x\in \bbP( V_\bbC) \st q(x)=0,\,\, q(x,\bar{x}) > 0\}
\end{equation}
and {\em the period map} $\Per\colon\Teich(M)\to\Perspace$, sending the point corresponding to a complex structure $I$ to $[\sigma_I]$.

Consider a connected component $\Teich^\circ(M)$ of the Teichm\"uller space, and let $\Per^\circ$ be the restriction of $\Per$ to $\Teich^\circ(M)$.  The subgroup $\MCG^\circ(M)\subset \MCG(M)$ that preserves the connected component $\Teich^\circ(M)$ is called {\em the monodromy group} of $\Teich^\circ(M)$.  Recall, see \cite{_Verbitsky:Torelli_}, that one can introduce an equivalence relation on the points of $\Teich(M)$, defining that $[I_1]\sim[I_2]$ when any open neighbourhood of $[I_1]$ intersects any open neighbourhood of $[I_2]$.  According to a theorem of Huybrechts, see \cite{_Huybrechts:basic_}, if $[I_1]\sim[I_2]$, then the complex manifolds $(M, I_1)$ and $(M,I_2)$ are bimeromorphic. Taking the quotient of $\Teich^\circ(M)$ by the above equivalence relation, one obtains a Hausdorff complex manifold $\Teich^\circ_\sim(M)$ (see \cite{_Verbitsky:Torelli_}).  Since the manifold $\Perspace$ is Hausdorff, the period map factors through $\Teich^\circ_\sim(M)$, and {\em the global Torelli theorem}, see \cite{_Verbitsky:Torelli_}, claims that $\Per^\circ_\sim\colon \Teich_\sim^\circ(M) \to \Perspace$ is an isomorphism of complex manifolds.

Next we give a more precise description of the fibres of $\Per^\circ$.  Following \cite{_AV:MBM_}, recall that for a hyperk\"ahler manifold $M$ one can define a collection of elements of $V_\bbZ = H^2(M,\bbZ)$, called {\em the MBM classes}, that are used to describe the K\"ahler cone of $M$.  The set of MBM classes may \textit{a priori} depend on the connected component $\Teich^\circ(M)$; we denote this set by $\MBM^\circ \subset V_\bbZ$.  Recall that $\MBM^\circ$ is invariant under the $\MCG^\circ(M)$-action and consists of a finite number of $\MCG^\circ(M)$-orbits.  The MBM classes can be characterized by the following property. Let $(M,I)$ be a deformation of $M$ such that its Picard group is generated by a primitive $x\in V_\Z$ with negative BBF square, and let $\eta\in H_2(M, \Q)$ be the BBF dual homology class.  Then $x$ is MBM if and only if $\eta$ is $\Q$-effective, that is, proportional to a homology class of a complex curve. Note that if this property holds for one such $I$, it holds for all $I$ such that $\Pic(M,I)=\langle x\rangle$ (see \cite{_AV:MBM_}).

For a complex structure $I\in \Comp(M)$, the Hodge decomposition on $V_\bbC = H^2(M,\bbC)$ depends only on the period $p = \Per(I)\in\Perspace$.  More precisely:
$$
V_\bbC = V^{2,0}_p\oplus V^{1,1}_p \oplus V^{0,2}_p,
$$
where $V^{2,0}_p$ is the subspace of $V_\bbC$ corresponding to the point $p$, $V^{0,2}_p$ is the complex conjugate of $V^{2,0}_p$ and $V^{1,1}_p$ is the $q$-orthogonal complement of $V^{2,0}_p\oplus V^{0,2}_p$.

Let $V^{1,1}_{p, \bbR} = V^{1,1}_p \cap H^2(M,\bbR)$. The restriction of the BBF form to $V^{1,1}_{p, \bbR}$ has signature $(1,b_2(M) - 3)$.  Let $\KK_I \subset V^{1,1}_{p,\bbR}$ be {\em the K\"ahler cone} of $(M,I)$, \textit{i.e.}~the collection of all K\"ahler classes on this complex manifold. Let $V_p^+ \subset V^{1,1}_{p,\bbR}$ be {\em the positive cone}, \textit{i.e.}~the connected component of the subset of $q$-positive classes that contains $\KK_I$. Recall that the MBM classes are $q$-negative and their orthogonal complements define a locally finite collection of walls in $V^+_p$ (see \cite{_SSV:rigid_currents_}).  Define $\MBM^{1,1}_p = \MBM^\circ\cap V^{1,1}_{p,\bbR}$.  As shown in \cite{_AV:MBM_}, the K\"ahler cone $\KK_I$ is a connected component of the set
\begin{equation}\label{eqn_chambers}
\CC^+_p = V^+_p\setminus \bigcup\limits_{x\in\MBM^{1,1}_p} x^\perp.
\end{equation}
We can now describe the fibres of the period map: for $p\in\Perspace$ we have $\Per^{-1}(p) = \pi_0(\CC^+_p)$. The connected components of $\CC^+_p$ are the K\"ahler cones of the hyperk\"ahler manifolds in the fibre $\Per^{-1}(p)$. These connected components are called {\em the K\"ahler chambers} of the positive cone.

For a point $p\in\Perspace$ let $N^{1,1}_p = V^{1,1}_p \cap V_\bbQ$ be the {\em N\'eron--Severi space}, \textit{i.e.}~the rational vector space spanned inside $V^{1,1}_p$ by the N\'eron--Severi group of a hyperk\"ahler manifold $(M, I)$ with $\Per(I) = p$.  The restriction of~$q$ to $N^{1,1}_p$ makes the latter a quadratic vector space.  If the restriction of~$q$ is non-degenerate, $N^{1,1}_p$ is either {\em hyperbolic}, \textit{i.e.}~of signature $(1, \rho -1)$, where $\rho$ is the Picard number of $(M,I)$, or {\em elliptic}, \textit{i.e.}~the restriction of~$q$ to $N^{1,1}_p$ is negative definite.  If the restriction of $q$ to the N\'eron--Severi space is degenerate, then it has one-dimensional kernel and is negative semidefinite. In this case $N^{1,1}_p$ is {\em parabolic}. By a well-known criterion of Huybrechts, see \cite{_Huybrechts:basic_}, the manifold $(M,I)$ is projective if and only if $N^{1,1}_p$ is hyperbolic.

\subsection{Semiample line bundles and Lagrangian fibrations}\label{sec_semiample}

Let $M$ be a compact complex manifold and $L\in\Pic(M)$. If $H^0(M,L)\neq 0$, then the canonical evaluation morphism $H^0(M,L)\otimes\OO_M\to L$ induces a meromorphic map $\varphi_L\colon M\dashrightarrow \bbP H^0(M,L)^*$.  We denote by $\im(\varphi_{L})$ the closure of its image.  We define {\em the Iitaka dimension} of $L$, or {\em the $L$-dimension} (see \cite{_Ueno_,_Lieberman_Sernesi_}) as
$$
\kappa(L) = \sup_{k > 0}\left\{ \dim \im \left(\varphi_{L^k}\right)\right\}.
$$
Here we use the convention that if $H^0(M,L^k) = 0$, then $\dim \im (\varphi_{L^k}) = -\infty$.  Recall that the {\em algebraic dimension} of $M$, denoted by $a(M)$, is the transcendence degree over $\bbC$ of the field of meromorphic functions on $M$. Since $\im(\varphi(L^k))$ is a projective variety, the rational functions on $\im(\varphi(L^k))$ form a subfield in the meromorphic functions on $M$, and we clearly have $\kappa(L) \le a(M)$.

Recall that the line bundle $L$ is called {\em nef}\, if $c_1(L)$ lies in the closure of the K\"ahler cone of $M$. The bundle~$L$ is called {\em semiample} if $L^k$ is generated by its global sections for some $k > 0$, \textit{i.e.}~$\varphi_{L^k}$ is a regular morphism.

Now assume that $M$ is a hyperk\"ahler manifold of dimension $2n$ and $I\in \Comp(M)$.  Assume that $L$ is nef, and let $\ell = c_1(L)$. Since $\KK_I$ is contained in the positive cone $V^+_p$, where $p = {\Per(I)}$, we always have $q(\ell)\ge 0$. If $q(\ell) > 0$, then $(M, I)$ is projective by the theorem of Huybrechts mentioned in the previous section.  Then a theorem of Kawamata \cite[Theorem 6.1]{_Kawamata:Pluricanonical_} implies that $L$ is semiample. Moreover, the Fujiki relations (\ref{eqn_Fujiki}) imply that $\int_M\ell^{2n} > 0$; therefore, the image of $\varphi_{L^k}$ has dimension $2n$, so that $\kappa(L) = \dim(M)$, \textit{i.e.}~$L$ is big.

We will be interested in the case when $L$ is semiample, $\ell = c_1(L)\neq 0$ and $q(\ell) = 0$.  In this case, after replacing $L$ by its power, we get a morphism $\pi\colon M\to B$ onto a projective variety $B$ such that $L = \pi^*\OO_B(1)$, where $\OO_B(1)$ is an ample line bundle on $B$. Passing to the normalization of $B$ and considering the Stein factorization, we may also assume that $B$ is normal and the fibres of $\pi$ are connected; \textit{i.e.}~$\pi$ is {\em a fibration} over a normal projective variety $B$.  Since in this case, by the Fujiki relations, $\int_M\ell^{2n} = 0$, we have $\dim(B) < \dim(M)$. Since $L$ is non-trivial and semiample, we also have $\dim(B) > 0$. Then the following fundamental result of Matsushita implies that $\pi$ is {\em a Lagrangian fibration}.

\begin{theorem}[\textit{cf.} Matsushita,  \cite{_Matsushita:fibre_,_Matsushita:equi_}]    \label{_Matsushita_fibra_Theorem_}
Let $\pi\colon M \to B$ be a surjective holomorphic map from a hyperk\"ahler manifold $M$ to a normal projective variety $B$, with $0< \dim B < \dim M$.  Then $\dim B = \frac{1}{2} \dim M$, and all fibres of $\pi$ are holomorphic Lagrangian subvarieties of $M$.
\end{theorem}

\begin{remark}
Theorem~\ref{_Matsushita_fibra_Theorem_} claims that all fibres of $\pi$, including the singular ones, are Lagrangian subvarieties. It means the following.  Assume that $Z\subset M$ is an irreducible component of one of the fibres of~$\pi$. Let $r\colon Z' \to Z\subset M$ be a resolution of singularities. Then Theorem~\ref{_Matsushita_fibra_Theorem_} claims that the dimension of $Z$ is $\frac{1}{2} \dim M$ and $r^*\sigma = 0$, where $\sigma$ is the holomorphic symplectic form on $M$. This claim is in \cite[Corollary~1]{_Matsushita:equi_}.
  \end{remark}

\begin{remark}
The base $B$ of a Lagrangian fibration is conjectured to be biholomorphic to $\C P^n$.  In \cite{_Matsushita:higher_} Matsushita has shown that $B$ is Fano and has the same rational cohomology as $\C P^n$ when it is smooth.  In \cite{_Hwang:base_} Hwang has shown that in this case $B$ is in fact biholomorphic to $\C P^n$.  In \cite{_Huybrechts_Xu_} Huybrechts and Xu have shown that $B$ is always biholomorphic to $\C P^2$ if it is normal and $\dim_\C M=4$ (see also \cite{_Ou:4-fold_,_Bogomolov_Kurnosov:4folds_}).
\end{remark}

Next recall that, by another result of Matsushita, the property of semiampleness of a line bundle is preserved under small deformations. More precisely, we have the following statement.

\begin{proposition}[\textit{cf.} Matsushita, \protect{\cite[Corollary 1.3]{_Matsushita:families_}}]  \label{prop_stability_sa}
  Let $\varphi\colon\MM\to T$ be a smooth family of hyperk\"ahler manifolds and $\LL\in\Pic(\MM)$ a line bundle.  Let $\MM_t = \varphi^{-1}(t)$ and $\LL_t = \LL|_{\MM_t}$ for $t\in T$. Assume that $\LL_{t_0}$ is semiample for some $t_0\in T$ and $q(c_1(\LL_{t_0})) = 0$. Then $\LL_t$ is semiample for all $t$ in some open neighbourhood of $t_0$.
  \end{proposition}

\subsection{C-symplectic structures and degenerate twistor deformations}

We recall the main results of \cite{_V:degenerate_,_SV:Moser_,_SV:Herm_} for later use. We assume that $\pi\colon M \to B$ is a Lagrangian fibration on a hyperk\"ahler manifold $M$ of dimension $2n$ over a normal projective base $B$. Denote by $\sigma \in \Lambda^2M\otimes \bbC$ the holomorphic symplectic form on $M$ and by $\eta\in \Lambda^{1,1}B$ a K\"ahler form on $B$. For $t\in \bbC$ let $\sigma_t = \sigma + t \pi^*\eta$.  One can check, see \cite{_V:degenerate_}, that $\sigma_t$ has the following properties:
\begin{enumerate}[(1)]
\item $d\sigma_t = 0$;
\item $\sigma_t^{n+1} = 0$;
\item $\sigma_t^n\wdg\bar{\sigma}_t^n$ is a volume form.
\end{enumerate}
We define {\em a C-symplectic form} as a complex-valued 2-form on a smooth manifold which satisfies these three assumptions.

A C-symplectic form defines a complex structure on $M$ as follows. Consider the map given by the contraction with $\sigma_t$:
$$
\iota_{\sigma_t}\colon TM\otimes \bbC \lra \Lambda^1M\otimes \bbC.
$$
One can check that $\ker(\iota_{\sigma_t})$ is a subbundle of rank $2n$; let us denote this subbundle by $T^{0,1}_t M$. Then one can prove that there exists a decomposition $$TM\otimes \bbC = T^{1,0}_t M \oplus T^{0,1}_t M,$$ where $T^{1,0}_t M$ is the complex conjugate of $T^{0,1}_t M$.  This defines an almost complex structure $I_t$ on $M$. Since the form $\sigma_t$ is closed, it is easy to see that $T^{0,1}_t M$ is closed under the Lie bracket; hence $I_t$ is integrable.  By construction, the form $\sigma_t$ is a holomorphic symplectic form on $(M, I_t)$.

The map $\pi$ remains holomorphic as a map from $(M,I_t)$ to $B$, where the complex structure on the base is unchanged.  It is also clear from the definition of $\sigma_t$ that the fibres of $\pi$ remain Lagrangian and retain the same complex structure. Therefore, the family of complex structures $I_t$ gives a deformation of the original Lagrangian fibration.

It is shown in \cite{_V:degenerate_} that there exists a smooth family $\varphi\colon\MM\to \bbC$ of complex manifolds such that $\MM_t \simeq (M,I_t)$ for the complex structures $I_t$ described above and any $t\in \bbC$.  We call the family $\MM$ {\em a degenerate twistor deformation} of $M$.  We also have the following result that will be important for us later.

\begin{theorem}[\textit{cf.} \protect{\cite[Theorem 1.1]{_SV:Herm_}}] \label{thm_kahler}
All fibres $\MM_t$ of the degenerate twistor family admit K\"ahler metrics.
\end{theorem}

Since the fibres $\MM_t$ also carry holomorphic symplectic forms, they are hyperk\"ahler by the Calabi--Yau theorem, \textit{i.e.}~$I_t\in \Comp(M)$. Therefore,  the family $\varphi$ defines an affine line in the Teichm\"uller space $\Teich(M)$.  We will call it {\em a degenerate twistor line}. Let us describe the image $F$ of this line in the period domain. We may assume that the form $\eta\in \Lambda^{1,1}B$ represents an integral cohomology class and let $\ell = [\pi^*\eta] \in V_\bbZ$. Let $W\subset V_\bbC$ be the subspace spanned by $[\sigma]$ and $\ell$.  Since $\Per(\sigma_t) = [\sigma + t\pi^*\eta]$, we see that $F$ is contained in the projective line $\bbP(W)\subset \bbP(V_\bbC)$ passing through the points $[\sigma]$ and $[\ell]$. In fact, it is clear from the definition (\ref{eqn_perspace}) of $\DD$ that~$F$ is the affine line $\bbP(W)\cap \DD = \bbP(W)\setminus\{[\ell]\}$.

\subsection{The semiample Teichm\"uller space}\label{sec_Teich_sa}

As before, we let $M$ be a hyperk\"ahler manifold.  As was explained above, the Teichm\"uller space of $M$ parametrizes the complex structures of hyperk\"ahler type on $M$ up to isotopy. Our goal is to introduce and study a space that parametrizes pairs $(I, L)$, where $I\in \Comp(M)$ and $L$ is a semiample line bundle on $(M, I)$, again up to isotopy.  We will call this space the semiample Teichm\"uller space (see~\Cref{defn_Teich_sa}).

We will always study the complex structures in some fixed connected component $\Teich^\circ(M)$ of the Teichm\"uller space.  It is clear that, within a connected component of the semiample Teichm\"uller space, the first Chern class of the bundle $L$ cannot change.  Therefore, we start by fixing a primitive non-zero cohomology class $\ell \in V_\bbZ = H^2(M,\bbZ)$ such that $q(\ell) = 0$. If $L\in\Pic(M, I)$ with $c_1(L) = \ell$, then $\ell$ is of Hodge type $(1,1)$ on $(M,I)$, so $\Per(I)$ is $q$-orthogonal to $\ell$.  We need to consider only those deformations of $I$ for which $\ell$ stays of Hodge type $(1,1)$, \textit{i.e.}~such that $L$ deforms together with $I$ as a holomorphic line bundle.  The periods of such deformations of $I$ are orthogonal to $\ell$. We therefore consider the intersection
$$
\widetilde\DD_\ell = \DD\cap \bbP(\ell^\perp).
$$
The domain $\widetilde\DD_\ell$ has two connected components; we will consider only one of them.  Below we explain this in more detail.

Note that the restriction of $q$ to $\ell^\perp\subset V_\bbR$ has one-dimensional kernel spanned by $\ell$. Let $W = \ell^\perp / \bbR \ell$ be the quotient with the induced quadratic form, which we also denote by $q$. The signature of $q$ on $W$ is $(2,\dim W - 2)$. It is well known (see \textit{e.g.}~\cite{_Verbitsky:Torelli_}) that the corresponding period domain
$$
\DD' = \{y\in \bbP(W\otimes \bbC)\st q(y) = 0,\, q(y,\bar{y})>0\}
$$
has two connected components.  Indeed, $\DD'$ is isomorphic to the homogeneous space
\[ \frac{O(2, \dim W-2)}{SO(2) \times O(\dim W-2)},\]
the group $O(2, \dim(W)-2)$ has four connected components, and the group $SO(2) \times O(\dim W-2)$ has two connected components.

Taking the quotient of $\ell^\perp\subset V_\bbC$ by $\bbC\ell$ induces a natural $\bbA^1$-fibration
$$
\widetilde\DD_\ell \lra \DD'; 
$$
therefore, the space $\widetilde\DD_\ell$ also has two connected components.  The following definition makes it possible to choose one of these two connected components unambiguously.

\begin{definition}\label{defn_perspace_l}
We define the {\em $\ell$-period domain} $\DD_\ell$ to be the connected component of $\widetilde\DD_\ell$ that satisfies the following condition: for some $[I]\in \Teich^\circ(M)$ with $p = \Per(I) \in \widetilde\DD_\ell$, the class $\ell$ lies in the closure of the positive cone $V_p^+$.
\end{definition}

We let $\Teich(M,\ell) = \Per^{-1}(\DD_\ell)$ be the corresponding divisor in
the Teichm\"uller space and $\Teich^\circ(M,\ell) = \Teich^\circ(M)\cap \Teich(M,\ell)$.

It is well known that the period map $\Teich^\circ(M)\to \Perspace$ is
generically one-to-one. The same is true for the map $\Per\colon
\Teich^\circ(M,\ell)\to \DD_\ell$.

\begin{proposition}\label{prop_generic}
Let $M$ be a hyperk\"ahler manifold, $\Teich^\circ(M,\ell)$ and $\DD_\ell$ as above. Consider a very general point $p\in \DD_\ell$. Then the preimage $\Per^{-1}(p)$ in $\Teich^\circ(M,\ell)$ is a single point. Moreover, the space $\Teich^\circ(M,\ell)$ is connected.
\end{proposition}

\begin{remark}\label{rem_generic}
  The proof given below actually establishes a stronger statement: for a very general point $p\in \DD_\ell$, the group of integral $(1,1)$-classes on $(M,I)$ is generated by $\ell$ for any $I\in \Per^{-1}(p)$.
  \end{remark}

\begin{proof}
We need to prove that, for a very general $p\in \DD_\ell$, the set $\CC^+_p$ is connected, \textit{i.e.}~that $\MBM^{1,1}_p$ is empty for such $p$. For a fixed $x\in \MBM^\circ$ we have $q(x)< 0$, and the subset $\DD_x = \DD\cap \bbP(x^\perp)$ is an irreducible divisor in $\DD$ (see \textit{e.g.}~\cite{_Verbitsky:Torelli_}).  Observe that $\DD_x\cap \DD_\ell$ is a proper subvariety of $\DD_\ell$, because $\langle x, l\rangle^\perp$ is of codimension~$2$ in $V_\bbR$. Therefore, for any point $p$ in
$$
\DD_\ell^\circ = \DD_\ell\setminus \bigcup\limits_{x\in \MBM^\circ} \DD_x, 
$$
we have $\MBM^{1,1}_p = \emptyset$. This proves the first claim.

For the second claim, recall that by~\Cref{defn_perspace_l} the $\ell$-period domain  $\DD_\ell$ is connected, and by the global Torelli theorem $\Teich^\circ(M,\ell)$ maps surjectively onto $\DD_\ell$. If $\Teich^\circ(M,\ell) = \TT_1\coprod\TT_2$, where $\TT_1$ and $\TT_2$ are non-empty open subsets, their images $U_1$ and $U_2$ are open subsets of $\DD_\ell$ (the period map is open) and $\DD_\ell = U_1\cup U_2$.  Since $\DD_\ell$ is connected, it follows that $U_1\cap U_2\neq\emptyset$, so $U_1\cap U_2\cap \DD_\ell^\circ\neq\emptyset$, and, by the first claim proven above, $\TT_1\cap \TT_2\neq\emptyset$, contradicting our assumption. This proves the second claim.
\end{proof}

Recall that $M$ is simply connected, so given $[I]\in \Teich(M,\ell)$ there exists a unique $L\in \Pic(M,I)$ with $c_1(L) = \ell$. Therefore, the semiample Teichm\"uller space may be defined as a subset of the usual one, as follows.

\begin{definition}\label{defn_Teich_sa}
In the above setting we define the {\em semiample Teichm\"uller space} as
\begin{eqnarray*}
\Teich_{\sa}(M,\ell) = \{[I]\in \Teich(M,\ell)&|& \mbox{for } L\in \Pic(M,I) \mbox{ with } c_1(L) = \ell, \\
&&  L \mbox{ is semiample}\}
\end{eqnarray*}
with the induced topology, and let
$$
\Per_{\sa}\colon\Teich_{\sa}(M,\ell)\lra\DD_{\ell}
$$
be the corresponding period map.
\end{definition}

\begin{proposition}\label{prop_Teich_open}
In the above setting $\Teich_{\sa}(M,\ell)$ is an open subset of\, $\Teich(M,\ell)$.
\end{proposition}

\begin{proof} We may assume that $\Teich_{\sa}(M,\ell)\neq \emptyset$; otherwise,  the claim clearly holds.
  The statement follows directly from \cite{_Matsushita:families_} (see also~\Cref{prop_stability_sa}).  For a point $[I]\in \Teich_{\sa}(M,\ell)$, we consider the universal deformation $\varphi\colon \MM\to U$ of $(M, I)$ whose base $U$ is an open neighbourhood of $[I]$ in $\Teich(M)$.  The intersection $U_\ell = U\cap \Teich(M,\ell)$ is an open neighbourhood of $[I]$ in $\Teich(M,\ell)$, and by~\Cref{prop_stability_sa} after shrinking $U_\ell$ we get $U_\ell \subset \Teich_{\sa}(M,\ell)$.
  \end{proof}

We will also use the following version of the Teichm\"uller space, which is \textit{a priori}  larger than the semiample Teichm\"uller space.

\begin{definition}\label{defn_Teich_nef}
We define the {\em nef Teichm\"uller space} as
\begin{eqnarray*}
\Teich_{\nef}(M,\ell) = \{[I]\in \Teich(M,\ell)&|& \mbox{for } L\in \Pic(M,I) \mbox{ with } c_1(L) = \ell,  \\
&&  L \mbox{ is nef}\}.
\end{eqnarray*}
\end{definition}

The stronger version of the SYZ conjecture (\Cref{_SYZ_strong_}) claims that $\Teich_{\nef}(M,\ell)=\Teich_{\sa}(M,\ell)$.

To study the fibres of $\Per_{\sa}$, it will be important to consider the collection of MBM classes orthogonal to the given class $\ell$. We fix a connected component $\Teich^\circ(M)$ of the Teichm\"uller space and the corresponding collection of the MBM classes $\MBM^\circ\subset V_\bbZ$.  Let $\MBM^\circ_\ell = \MBM^\circ \cap \ell^\perp$ and $\MCG^\circ_\ell(M)$ be the stabilizer of $\ell$ in the monodromy group $\MCG^\circ(M)$.

Given a point $p\in \DD_\ell$, we define $\MBM_{\ell, p}^{1,1} = \MBM_\ell \cap V^{1,1}_p$.  The classes in $\MBM_{\ell, p}^{1,1}$ define a wall-and-chamber decomposition of the positive cone:
\begin{equation}\label{eqn_chamber_stable}
\CC^+_{\ell,p} = V^{+}_p \setminus \bigcup\limits_{x\in\MBM^{1,1}_{\ell,p}} x^\perp.
\end{equation}
We will call the connected components of $\CC^+_{\ell,p}$ the {\em $\ell$-stable K\"ahler chambers} of the positive cone.

\begin{figure}[h]
\centering
\begin{minipage}{.45\textwidth}
  \centering
  \includegraphics[width=\linewidth]{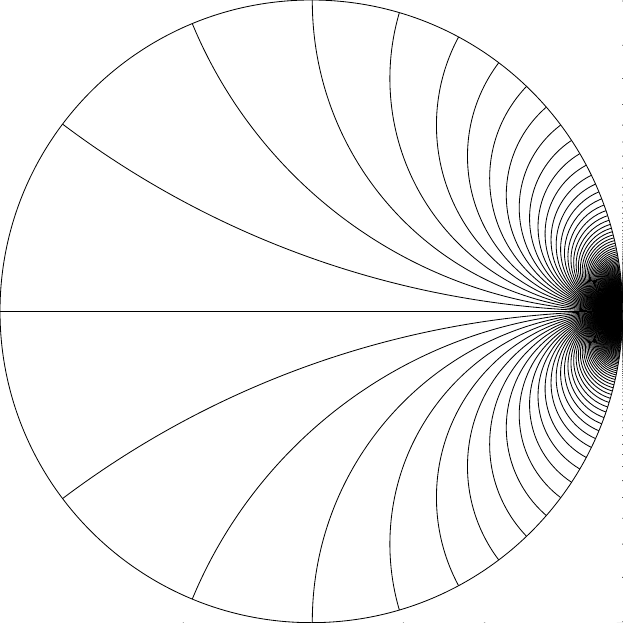}
  \captionof{figure}{Stable K\"ahler chambers}
  \label{fig_1}
\end{minipage}\ \ \ \ \ 
\begin{minipage}{.45\textwidth}
  \centering
  \includegraphics[width=\linewidth]{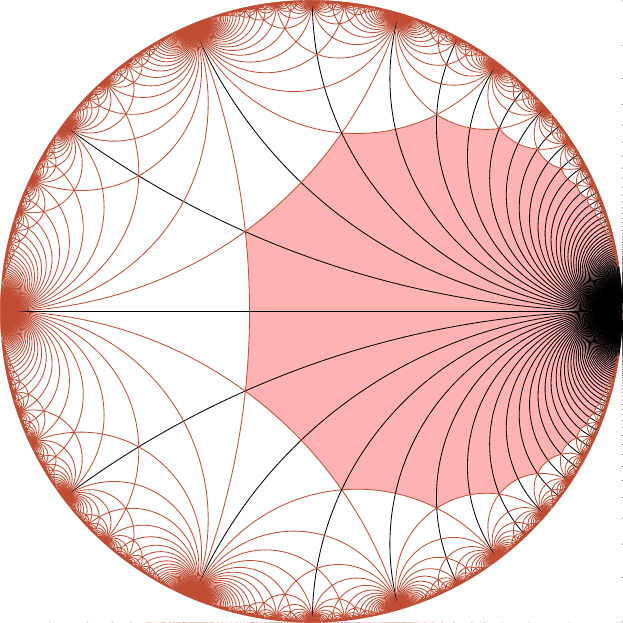}
  \captionof{figure}{K\"ahler chambers}
  \label{fig_2}
\end{minipage}
\end{figure}

The difference between the K\"ahler chambers and the $\ell$-stable K\"ahler chambers is illustrated on Figures~\ref{fig_1} and~\ref{fig_2}. The projectivisation of $V^+_p$ is a hyperbolic space, and the orthogonal complements of the MBM classes define a collection of walls in that space. The $\ell$-stable K\"ahler chambers are cut out by the walls that pass through the point $[\ell]$, as shown on Figure~\ref{fig_1}, where the point $[\ell]$ is the rightmost point of the boundary, where the walls accumulate. Each stable K\"ahler chamber is further cut into K\"ahler chambers by the walls that do not pass through $[\ell]$, as shown on Figure~\ref{fig_2}. Every $\ell$-stable K\"ahler chamber contains a unique K\"ahler chamber that has $[\ell]$ in its closure. These K\"ahler chambers are shown in light red on Figure~\ref{fig_2}; they correspond to the points $[I]\in \Teich^\circ(M,\ell)$ such that $\ell$ is nef on $(M,I)$.

As we will prove in~\Cref{prop_fibered}, the $\ell$-stable K\"ahler chambers, \textit{i.e.}~the connected components of $\CC^+_{\ell,p}$, are in bijection with the degenerate twistor lines in the Teichm\"uller space passing over the point $p$.  We have the following observation.

\begin{proposition}\label{prop_inf_or_none}
The set of $\ell$-stable K\"ahler chambers is infinite when $\MBM^\circ_\ell$ is non-empty and consists of one element when $\MBM^\circ_\ell$ is empty.
\end{proposition}

\begin{proof} 
The number of connected components of $\CC^+_{\ell,p}$ is infinite whenever the number of walls is infinite. Therefore, it is enough to show that the set $\MBM^\circ_\ell$ is either infinite or empty. We will construct an infinite group that acts on this set with all orbits infinite.

Denote by $\Gamma$ the image of $\MCG^\circ(M)$ in the orthogonal group $O(V_\bbQ,q)$.  Let $G$ be the stabilizer of $\ell$ in $O(V_\bbQ,q)$.  Fix a decomposition of $V_\bbQ$ into an orthogonal sum of two $\bbQ$-vector spaces $V_\bbQ = U\oplus W$, where~$U$ is a 2-dimensional subspace of signature (1,1) that contains $\ell$ and $W = U^\perp$.  The group $G$ acts on $\ell^\perp = W\oplus \bbQ\ell$, so it also acts on $W$. Let $K$ be the kernel of the homomorphism $G\to O(W,q)$. The group~$K$ is abelian, and if we identify the Lie algebra of $O(V_\bbQ,q)$ with $\Lambda^2 V_\bbQ$, then the Lie algebra of $G$ is $\Lambda^2(\ell^\perp) \simeq \ell \wdg W \oplus \Lambda^2 W$ and the Lie algebra of $K$ is $\ell\wdg W$.  The elements of $K$ that preserve the integral structure $V_\bbZ$ form a lattice in $K$ that we denote by $K_\bbZ$. It is a free abelian group of rank $\dim W$.  Let us denote by $\Gamma_\ell$ the intersection of~$\Gamma$ with $K_\bbZ$. It is known that $\Gamma$ is a subgroup of finite index in $O(V_\bbZ,q)$; see \textit{e.g.}~\cite{_Verbitsky:Torelli_}. Therefore, $\Gamma_\ell$ is of finite index in $K_\bbZ$; \textit{i.e.}~$\Gamma_\ell$ is also a free abelian group of rank $\dim W$. The group $\Gamma_\ell$ acts on $\ell^\perp$ by translations in the direction of $\ell$.  More precisely, given $y\in W$, the element $\ell\wdg y \in \ell\wdg W \simeq \mathrm{Lie}(K)$ acts on $\ell^\perp$ as follows:
$$x\longmapsto q(x,y)\ell - q(x,\ell)y = q(x,y)\ell.$$ Exponentiating, we obtain the following formula for the action of $\gamma = e^{\ell\wdg y} \in K$ on $\ell^\perp$:
\begin{equation}\label{_para_action_Equation_}
\gamma\colon x\longmapsto x+ q(x, y) \ell.
\end{equation}
The subgroup $\Gamma_\ell$ is identified with a lattice $R \simeq \Z^{\dim W}\subset W$, so that $\gamma = e^{\ell\wdg y} \in \Gamma_\ell$ if and only if $y \in R$.

Given a point $p\in \DD_\ell$, the subspace $V'_\bbC$ spanned by $p$, $\bar{p}$ and $\ell$ is contained in $\ell^\perp$. Observe that $\Gamma_\ell$ is contained in $K$ by construction; therefore, it preserves the subspace $V'_\bbC$ and its orthogonal complement. Being a subgroup of $\Gamma$, it also preserves the set of MBM classes. Hence $\Gamma_\ell$ preserves $\MBM^\circ_\ell$.  Clearly, the action \eqref{_para_action_Equation_} has an infinite orbit whenever $x\notin W^\bot=\Q\ell$.  Therefore, any element of $\MBM^\circ_\ell$ has an infinite $\Gamma_\ell$-orbit.  This completes the proof.
\end{proof}

\section{The main results}

In this section we state and prove the main results describing the structure of the semiample Teichm\"uller space introduced above.  We let $M$ be a hyperk\"ahler manifold and $\Teich^\circ(M)$ a connected component of the Teichm\"uller space. We let $\Teich^\circ_{\sa}(M) = \Teich_{\sa}(M,\ell) \cap \Teich^\circ(M)$ and $\Per^\circ_{\sa}\colon \Teich^\circ_{\sa}(M)\to \DD_\ell$ be  the restriction of the period map. Analogously, let $\Teich^\circ_{\nef}(M) = \Teich_{\nef}(M,\ell) \cap \Teich^\circ(M)$.

\subsection{Algebraic dimension of the manifolds with periods in $\boldsymbol{\DD_\ell}$}

Recall that we denote by $a(M)$ the algebraic dimension of a compact complex manifold $M$. We have $a(M) \leq \dim_\C M$, see \cite{_Ueno_}, and the equality is realized if and only if $M$ is Moishezon, \textit{i.e.}~bimeromorphic to a projective manifold.  In general, the algebraic dimension is hard to control. For example it is not semicontinuous in families; see \cite{_FP:non-semicontinuity_}.  It is conjectured (see \cite{_COP:non-alge_}) that, for a non-projective $2n$-dimensional hyperk\"ahler manifold $M$, we have either $a(M) = 0$ (when the N\'eron--Severi space of $M$ is elliptic) or $a(M) = n$ (when the N\'eron--Severi space of $M$ is parabolic). The following statement partially confirms this expectation.

\begin{proposition}\label{prop_algdim}
Assume that $\Teich^\circ_{\sa}(M)$ is non-empty.  Then for any point $[I]\in \Teich^\circ(M,\ell)$, we have $a(M, I)= n$ when $(M,I)$ is non-projective, and $a(M, I)= 2n$ when $(M,I)$ is projective.
\end{proposition}

\begin{proof}
  Recall that in Section~\ref{sec_semiample} we  introduced the notation $\kappa(L)$ for the Iitaka dimension of a line bundle~$L$.  For a point $[I] \in \Teich^\circ(M,\ell)$, denote by $L_I$ the line bundle on $(M, I)$ with $c_1(L_I) = \ell$. Define a subset $\TT\subset \Teich^\circ(M,\ell)$ as follows: $[I]\in \TT$ if and only if there exists an open neighbourhood $U\subset \Teich^\circ(M,\ell)$ containing $[I]$ such that for any $[I']\in U$ we have $\kappa(L_{I'}) \ge n$. By our assumptions and~\Cref{prop_Teich_open}, $\TT$ is a non-empty open subset of $\Teich^\circ(M,\ell)$.  Assume that $[J]$ is a point in the closure of $\TT$, and let $\varphi\colon\MM\to U_0$ be the universal deformation of $(M, J, L_{J})$, where $U_0$ is a connected open neighbourhood of $[J]$.  The subset $\TT\cap U_0$ is non-empty and open in $U_0$, and for $[I] \in \TT\cap U_0$ we have $\kappa(L_I) \ge n$. Applying \cite[Section~1, first theorem]{_Lieberman_Sernesi_}  to the family $\MM$, we see that there exist an integer $k$ and a subset $W\subset U_0$ such that $W$ is the complement of the union of a countable number of proper closed subvarieties of $U_0$, $\kappa(L_I) \ge k$ for all $[I] \in U_0$ and $\kappa(L_I) = k$ for $[I] \in W$. Since $\TT\cap U_0$ is non-empty and open in $U_0$, we have $W \cap \TT \neq \emptyset$. We conclude that $k \ge n$ and hence $\kappa(L_I) \ge n$ for all $[I]\in U_0$.  It follows that $\TT$ is closed in $\Teich^\circ(M,\ell)$, and since the latter space is connected by~\Cref{prop_generic}, we have $\TT = \Teich^\circ(M,\ell)$. Next we use the inequality $a(M, I)\ge \kappa(L_I)$ that was recalled above to prove the first claim. The second claim now follows from \cite[Theorem 3.6]{_COP:non-alge_}.
  \end{proof}

\begin{corollary}\label{cor_sa_dense}
Assume that $\Teich^\circ_{\sa}(M)$ is non-empty. Let $[I] \in \Teich^\circ(M,\ell)$ be a point such that $(M, I)$ is non-projective. Then $[I]\in \Teich^\circ_{\sa}(M)$.
\end{corollary}

\begin{proof}
  This is \cite[Theorem 3.4]{_COP:non-alge_}. We give an alternative proof, since our argument will be useful below.  Let $p = \Per(I)\in\DD_\ell$.  By our assumption the N\'eron--Severi lattice of $(M,I)$ is parabolic; therefore, all elements of $N^{1,1}_p$ are orthogonal to $\ell$.  It follows that $\MBM^{1,1}_{p} = \MBM^{1,1}_{\ell,p}$ and $\CC^+_p = \CC^+_{\ell,p}$; \textit{i.e.}~the K\"ahler chambers of the positive cone coincide with the $\ell$-stable K\"ahler chambers. The vector $\ell$ lies in the closure of any $\ell$-stable K\"ahler chamber; hence $\ell$ is a nef class on $(M, I)$. The conclusion follows from~\Cref{prop_algdim} and \cite[Theorem 3.7]{_COP:non-alge_}.
  \end{proof}

\subsection{$\boldsymbol{\bbA^1}$-fibration on the semiample Teichm\"uller space}

Recall from~\Cref{defn_perspace_l} of the period domain $\DD_\ell$ that there exists an $\bbA^1$-fibration
\begin{equation}\label{eqn_alpha}
\alpha\colon\DD_\ell\lra \DD'
\end{equation}
induced by the projection from $\ell^\perp\subset V_\bbC$ to the quotient $W_\bbC = \ell^\perp/\bbC \ell$.  If $p\in \DD_\ell$, then the fibre of $\alpha$ passing through $p$ is the intersection of $\DD_\ell$ with the projective line in $\bbP(V_\bbC)$ passing through the points $p$ and $[\ell]$. So the fibres of $\alpha$ represent the degenerate twistor families in the period domain.

Recall that the period domain $\DD_\ell\subset {\mathbb P} (V_\C)$ is defined as the set of all points in the quadric $q(x)=0$ which satisfy $q(x,\bar x) >0$. The fibre $\alpha^{-1}(\alpha(p))$ is an affine line contained in the intersection of $\DD_\ell$ and ${\mathbb P}( W_\C)$, where $W_\C\subset V_\C$ is the 3-dimensional space generated by $p, \bar p$ and $\ell$. Since $q$ is positive semidefinite on $\Re(W_\C)$, the intersection of $\DD_\ell$ and ${\mathbb P} (W_\C)$ is the union of two complex lines $R$, $\bar R$ minus the point $\ell$ where they intersect.  Both complements $F=R \backslash \{\ell\}$ and $\bar F=\bar R \backslash \{\ell\}$ are fibres of $\alpha$, with $F=\alpha^{-1}(\alpha(p))$ and $\bar F=\alpha^{-1}(\alpha(\bar p))$.

Let $F\simeq \bbA^1$ be a fibre of $\alpha$. For any two points $p_1, p_2 \in F$, we define a canonical identification between the stable K\"ahler chambers (\ref{eqn_chamber_stable})
\begin{equation}\label{_identi_chambers_Equation_}
\chi_{p_1 p_2}\colon \pi_0\left(\CC^+_{\ell,p_1}\right) \overset{\lowsim}\lra \pi_0\left(\CC^+_{\ell,p_2}\right)
\end{equation}
as follows. 

Given $x\in \MBM^\circ_\ell$, observe that $p_1\in x^\perp$ if and only if $p_2\in x^\perp$: this follows from the fact that $\ell\in x^\perp$ and the three points $p_1, p_2$ and $[\ell]$ lie on the same line in $\bbP(V_\bbC)$, so $x$ is orthogonal to any two of them if and only if $x$ is orthogonal to all three.  We conclude that $\MBM^\circ_{\ell,p_1} = \MBM^\circ_{\ell,p_2}$.  Let us denote this set by $\MBM^\circ_{\ell, F}$.

Let $\CC_1$ be a connected component of $\CC^+_{\ell,p_1}$. Any element $x\in \MBM^\circ_{\ell,F}$ determines two open half-spaces in $V_\bbR$: one where the linear form $x\cnv q$ is positive, and the other where it is negative.  Let $\HH_x \subset V_\bbR$ be the half-space containing $\CC_1$. Then $\CC_1$ is the intersection of $V^+_{p_1}$ and the subset
$$
\HH = \left(\bigcap\limits_{x\in \MBM^\circ_{\ell,F}}\HH_x\right) \subset V_\bbR.
$$

Now, the intersection $V^+_{p_2}\cap \HH$ is a connected component $\CC_2$ of $\CC^+_{\ell,p_2}$, and we define $\chi_{p_1 p_2}(\CC_1) = \CC_2$.  This gives the identification $\chi_{p_1 p_2}$ as claimed above.  Next, let us observe that the constructed identification preserves the K\"ahler cones of the manifolds in the degenerate twistor families. More precisely, we have the following statement.

\begin{proposition}\label{prop_chambers}
Given a degenerate twistor family $\varphi\colon \MM\to \bbC$, let $t_1, t_2\in\bbC$ be two points and $p_1, p_2\in \DD_\ell$ be the periods of $\MM_{t_1}$ and $\MM_{t_2}$.  Let $\CC\in \pi_0(\CC^+_{\ell,p_1})$ be the stable K\"ahler chamber containing the K\"ahler cone $\KK_{t_1}$ of $\MM_{t_1}$. Then the K\"ahler cone $\KK_{t_2}$ of $\MM_{t_2}$ is contained in $\chi_{p_1 p_2}(\CC)$, where $\chi_{p_1 p_2}$ is the bijective correspondence \eqref{_identi_chambers_Equation_} between the sets of stable K\"ahler chambers constructed above.
\end{proposition}

\begin{proof}
Let $F$ be the affine line in $\DD_\ell$ corresponding to the family $\varphi$.  We use the notation introduced above. It is enough to check that for any $x\in \MBM^\circ_{\ell,F}$ the K\"ahler cone $\KK_{t_1}$ is contained in $\HH_x$ if and only if $\KK_{t_2}$ is contained in $\HH_x$.

Let $\varkappa\colon \KK\to \bbC$ be the relative K\"ahler cone of the family $\MM$, \textit{i.e.}~the subset of $V_\bbR\times\bbC$ such that $\KK_t$ is the K\"ahler cone of $\MM_t$ for all $t\in \bbC$. It is well known that $\KK$ is an open subset of the bundle $\VV^{1,1}_\bbR$ of real classes of type $(1,1)$.  Given $x\in \MBM^\circ_{\ell,F}$, let $\KK^+ = \{(y,t) \in \KK\st q(y,x) > 0\}$ and $\KK^- = \{(y,t) \in \KK\st q(y,x) < 0\}$.  Then $\varkappa(\KK^+)$ and $\varkappa(\KK^-)$ are disjoint open subsets that cover $\bbC$, so one of them is empty. We conclude that~$\KK$ is contained either in $\HH_x\times\bbC$ or in its complement. Therefore, $\KK_{t_1}$ is contained in $\HH_x$ if and only if $\KK_{t_2}$ is contained in $\HH_x$, which completes the proof.
\end{proof}

Note that $\Teich^\circ_{\sa}(M)$ is fibred by the degenerate twistor lines. Indeed, for any point $I\in \Teich^\circ_{\sa}(M)$, the manifold $(M,I)$ is equipped with a Lagrangian fibration $\pi$ associated with the morphism $(M,I) \to {\mathbb P}(H^0(M, L_I^k)^*)$, where $L_I$ is the line bundle with $c_1(L_I)=\ell$. Semiampleness of $L_I$ is equivalent to existence of such a fibration; however, the degenerate twistor deformation preserves the Lagrangian fibration, so semiampleness of $L_I$ is retained.

In the rest of this subsection, we are going to enumerate the degenerate twistor lines in $\Teich^\circ_{\sa}(M)$ mapped to the same line in the period space; see~\Cref{prop_fibered}.

\begin{remark}
Let $\Tw(M) \to \C P^1$ be the twistor family associated with a hyperk\"ahler structure. It is well known (see \textit{e.g.}~\cite{_Verbitsky:alge_family_}) that a very general fibre $(M,I)$ in this family is non-algebraic, and, moreover, the BBF form is negative definite on its N\'eron--Severi lattice $H^{1,1}(M, I) \cap H^2(M,\Q)$. A similar statement is true for the degenerate twistor deformation: a very general fibre is non-algebraic, and the BBF form is negative semidefinite on $H^{1,1}(M, I) \cap H^2(M,\Q)$.
\end{remark}

\begin{lemma}\label{lem_generic}
Let $F\subset \DD_\ell$ be a fibre of $\alpha$. For a very general $p\in F$, the N\'eron--Severi space $N^{1,1}_p$ is parabolic.
\end{lemma}

\begin{proof}
Note that $N^{1,1}_p$ is parabolic if and only if it does not contain $q$-positive elements of $V_\bbZ$. Given $x\in V_\bbZ$ with $q(x) > 0$, assume that $F\subset \bbP(x^\perp)$. Since $[\ell]\in \overline{F}\subset \bbP(V_\bbC)$, we then have $[\ell] \in \bbP(x^\perp)$, \textit{i.e.}~$q(x,\ell) = 0$. This equality implies that $\bbP(\langle x, \ell \rangle^\perp)\cap \DD = \emptyset$, which gives a contradiction. We conclude that $F\cap \bbP(x^\perp)$ either is empty or consists of one point. Therefore, for a very general $p\in F$, the space $N^{1,1}_p$ does not contain $q$-positive elements, which completes the proof.
\end{proof}

\begin{proposition}\label{prop_fibered}
Assume that $\Teich^\circ_{\sa}(M)$ is non-empty. Let $p\in\DD_\ell$ be a point in the image of $\Teich^\circ_{\sa}(M)$ under the period map, and let $F\simeq \bbA^1$ be the fibre of $\alpha$ passing through $p$. Then the preimage of $F$ in $\Teich^\circ_{\sa}(M)$ is a disjoint union of affine lines mapping surjectively onto $F$ and indexed by $\pi_0(\CC^+_{\ell,p})$. For every connected component $\CC\in \pi_0(\CC^+_{\ell,p})$, there exists a unique $[I]\in \Teich^\circ_{\sa}(M)$ such that $\Per(I) = p$ and the K\"ahler cone of $(M,I)$ is contained in $\CC$.\looseness=-1
\end{proposition}

\begin{proof}
By~\Cref{lem_generic} we may find a point $p'\in F$ with parabolic $N^{1,1}_{p'}$, and by the global Torelli theorem, the preimage of $p'$ in $\Teich^\circ(M,\ell)$ consists of the points indexed by $\pi_0(\CC^+_{p'})$. But since $N^{1,1}_{p'}$ is parabolic, $\CC^+_{p'} = \CC^+_{\ell,p'}$ (see the proof of~\Cref{cor_sa_dense}). For any point $[I']$ in the preimage of $p'$, the manifold $(M,I')$ is not projective, so we have $[I']\in \Teich^\circ_{\sa}(M)$ by~\Cref{cor_sa_dense}. Through each such point $[I']$ passes a unique degenerate twistor line that maps surjectively onto $F$. These lines are indexed by $\pi_0(\CC^+_{\ell,p'})$, which is identified with $\pi_0(\CC^+_{\ell,p})$ by the map $\chi_{p'p}$ constructed above. The last claim follows from~\Cref{prop_chambers}.
\end{proof}

\subsection{A global Torelli theorem for the semiample Teichm\"uller space}

It is convenient to reformulate and summarize the results of the above discussion as a version of the global Torelli theorem for semiample Teichm\"uller spaces.

\begin{theorem}\label{thm_Torelli_sa}
Assume that $\Teich_{\sa}^\circ(M)$ is non-empty. Then it has the following properties:
\begin{enumerate}
\item\label{thm_T_s-1} Two points $[I_1], [I_2]\in \Teich_{\sa}^\circ(M)$ are non-separated in $\Teich_{\sa}^\circ(M)$
if and only if they are non-separated in $\Teich^\circ(M)$. 
\item\label{thm_T_s-2} The map $\Per^\circ_{\sa}$ induces an isomorphism of the Hausdorff reduction
of\, $\Teich_{\sa}^\circ(M)$ and $\DD_\ell$. In particular, $\Per^\circ_{\sa}$ is surjective. 
\item\label{thm_T_s-3} $\Teich_{\sa}^\circ(M)$ is connected.
\item\label{thm_T_s-4} For $p\in \DD_\ell$ the preimage of $p$ under $\Per^\circ_{\sa}$ is naturally identified
with $\pi_0(\CC^+_{\ell,p})$. In particular, a very general $p\in \DD_\ell$ has a unique preimage
in $\Teich^\circ_{\sa}(M)$.
\item\label{thm_T_s-5} $\Teich_{\sa}^\circ(M) = \Teich_{\nef}^\circ(M)$.
\end{enumerate}
\end{theorem}

\begin{proof}
{\em Part}~\eqref{thm_T_s-1}. If $[I_1]$ and $[I_2]$ are non-separated in $\Teich_{\sa}^\circ(M)$, then they are clearly non-separated in $\Teich^\circ(M)$. Conversely, let $[I_1]$ and $[I_2]$ be non-separated in $\Teich^\circ(M)$. Then they have the same period $p\in\DD_\ell$. Let $U_1$ and $U_2$ be open neighbourhoods of $[I_1]$ and $[I_2]$ in $\Teich^\circ(M)$ and $W = \Per(U_1)\cap \Per(U_2)\cap \DD_\ell$.  Then~$W$ is an open neighbourhood of $p$ in $\DD_\ell$, because the period map an open morphism. By~\Cref{prop_generic} and~\Cref{rem_generic}, for a very general point $p'\in W$, its preimage in $\Teich^\circ(M,\ell)$ consists of a single point $[I]$ such that $(M, I)$ is non-projective. By the uniqueness of the preimage of $p'$, we must have $[I]\in U_1\cap U_2$.  The point $[I]$ is contained in $\Teich_{\sa}^\circ(M)$ by~\Cref{cor_sa_dense}.  So $U_1\cap U_2\cap \Teich^\circ_{\sa}(M)$ is non-empty; hence $[I_1]$ and $[I_2]$ are non-separated in $\Teich_{\sa}^\circ(M)$.

{\em Part}~\eqref{thm_T_s-2}.  Let us first prove the surjectivity of $\Per^\circ_{\sa}$. For $p\in \DD_\ell$ let $F\subset \DD_\ell$ be the fibre of $\alpha$ passing through $p$; see (\ref{eqn_alpha}). By~\Cref{lem_generic} there exists a $p'\in F$ with parabolic $N^{1,1}_{p'}$, and by the global Torelli theorem there exits an $[I']\in \Teich^\circ(M,\ell)$ with $\Per(I') = p'$. Since $N^{1,1}_{p'}$ is parabolic, the manifold $(M,I')$ is non-projective, and by~\Cref{cor_sa_dense}, $[I']\in \Teich^\circ_{\sa}(M)$. By~\Cref{prop_fibered}, $F$ is the image of a degenerate twistor line in $\Teich^\circ_{\sa}(M)$, so $\Per^\circ_{\sa}$ is surjective.

It follows from~\eqref{thm_T_s-1} that the Hausdorff reduction of $\Teich^\circ_{\sa}(M)$ is an open subset of the Hausdorff reduction of $\Teich^\circ(M,\ell)$. The latter is isomorphic to $\DD_\ell$ via the period map, and since $\Per^\circ_{\sa}$ is surjective, the Hausdorff reduction of $\Teich^\circ_{\sa}(M)$ is also isomorphic to $\DD_\ell$.

{\em Part}~\eqref{thm_T_s-3}. Assume that $\Teich^\circ_{\sa}(M) = \TT_1\coprod \TT_2$, where $\TT_1$ and $\TT_2$ are open subsets. Then the sets $U_i = \Per^\circ_{\sa}(\TT_i)$ are non-empty open subsets of $\DD_\ell$, and by~\eqref{thm_T_s-2} they cover $\DD_\ell$. Since $\DD_\ell$ is connected, $U_1\cap U_2$ is non-empty, and therefore by~\Cref{prop_generic} it contains points with unique preimage in $\Teich^\circ(M,\ell)$.  It follows that $\TT_1$ intersects $\TT_2$, contradicting our assumption. Therefore, $\Teich^\circ_{\sa}(M)$ is connected.

{\em Part}~\eqref{thm_T_s-4}. 
This follows directly from~\eqref{thm_T_s-2} and Propositions~\ref{prop_fibered} and~\ref{prop_generic}.

{\em Part}~\eqref{thm_T_s-5}.  Let $[I]\in \Teich^\circ_{\nef}(M,\ell)$ and $p = \Per(I)$. Recall that the K\"ahler cone $\KK_I$ lies in some $\ell$-stable K\"ahler chamber $\CC\in \pi_0(\CC^+_{\ell,p})$.  By~\Cref{prop_fibered} there exists a unique $[I']\in \Teich^\circ_{\sa}$ with $\Per(I') = p$ and $\KK_{I'}\subset \CC$. Now observe that $\ell$ is nef both for $I$ and $I'$, so $\ell$ lies in the closure of both $\KK_I$ and $\KK_{I'}$. Both of these cones are subchambers of $\CC$ cut out by the MBM classes not orthogonal to $\ell$.  For every such MBM class $x$, one and only one of the half-spaces $\{y\in V_\R \mid q(x, y)>0\}$ and $\{y\in V_\R \mid q(x, y)<0\}$ contains $\ell$.  Since the closure of~$\KK_I$ and the closure of $\KK_{I'}$ contain $\ell$, this implies that $\KK_I = \KK_{I'}$. By the usual global Torelli theorem $[I] = [I']$; this finishes the proof.
\end{proof}

\subsection{The hyperk\"ahler SYZ in families}

\begin{theorem}\label{thm_abundance}
Let $M$ be a hyperk\"ahler manifold and $L\in \Pic(M)$ a nef line bundle with $q(c_1(L)) = 0$. Assume that the pair $(M,L)$ admits a deformation $(M',L')$ such that $L'$ is semiample. Then $L$ is semiample.
\end{theorem}

\begin{proof}
We let $\ell = c_1(L)\in H^2(M,\bbZ)$.  By our assumptions there exist a smooth family of hyperk\"ahler manifolds $\varphi\colon \MM\to T$ over a connected base $T$, a line bundle $\LL\in\Pic(\MM)$ and two points $t_1, t_2 \in T$ such that $(M,L) \simeq (\MM_{t_1},\LL_{t_1})$ and $(M',L') \simeq (\MM_{t_2},\LL_{t_2})$.  Passing to the universal cover, we may assume that $T$ is simply connected.  This gives a morphism $T\to \Teich^\circ(M,\ell)$ to some connected component of the Teichm\"uller space. Since $L'$ is semiample, we see that $\Teich^\circ_{\sa}(M,\ell)$ is non-empty, and we may apply~\Cref{thm_Torelli_sa}. Since $L$ is nef, $M$ lies in $\Teich_\nef^\circ(M,\ell)$.  By item~\eqref{thm_T_s-5} in~\Cref{thm_Torelli_sa}, $\Teich_\nef^\circ(M,\ell) = \Teich_\sa^\circ(M,\ell)$, so $L$ is semiample. \end{proof}


\end{document}